\documentclass[12pt, reqno]{amsart}
\allowdisplaybreaks
\begin{document}
\setcounter{page}{1}

\title[\hfilneg \hfil Ostrowski and Cebysev inequality  ]
{Some Ostrowski and Cebysev type inequality in three variables  on Time Scales}

\author[Deepak B. Pachpatte\hfil \hfilneg]
{Deepak B. Pachpatte}

\address{Deepak B. Pachpatte \newline
 Department of Mathematics,
 Dr. Babasaheb Ambedkar Marathwada University, Aurangabad,
 Maharashtra 431004, India}
\email{pachpatte@gmail.com}

\subjclass[2010]{26E70, 34N05, 26D10}
\keywords{Ostrowski inequality, three variables, time scales.}

\begin{abstract}
           The main objective of this paper is to study some Ostrowski  and Cebysev  type inequalities in three variables  on Time Scales.

\end{abstract}
\maketitle

\section{Introduction}
 In 1938 A.M. Ostrowski \cite{Ost} proved the inequality: 
Let $f:[a,b] \to \mathbb{R}$ be continuous on [a,b] and differentiable on $(a,b)$ with derivative $f':[a,b] \to \mathbb{R}$ being bounded on $(a,b)$, that is $\left\| {f'} \right\|_\infty   = \sup _{x \in \left( {a,b} \right)} \left| {f'\left( x \right)} \right| < \infty 
$. Then
\[
\left| {f\left( x \right) - \frac{1}{{b - a}}\int\limits_a^b {f\left( t \right)dt} } \right| \le \left[ {\frac{1}{4} + \frac{{\left( {x - \frac{{a + b}}{2}} \right)^2 }}{{\left( {b - a} \right)^2 }}} \right]\left( {b - a} \right)\left\| {f'} \right\|_\infty,  
\]
for all $x \in [a,b]$. In this the inequality gives an upper bound for the approximation of the integral average $\frac{1}{{b - a}}\int\limits_a^b {f\left( t \right)dt} $ by the value $f(x)$ at the point $x \in [a,b]$.

In 1882, P.L. Cebysev \cite{Ceb} proved the following inequality:

\[
\begin{array}{l}
 \left| {\frac{1}{{b - a}}\int\limits_a^b {f\left( x \right)g\left( x \right)dx}  - \left( {\frac{1}{{b - a}}\int\limits_a^b {f\left( x \right)dx} } \right)\left( {\frac{1}{{b - a}}\int\limits_a^b {g\left( x \right)dx} } \right)} \right| \\ 
  \le \frac{1}{{12}}\left( {b - a} \right)^2 \left\| {f'} \right\|_\infty  \left\| {g'} \right\|_\infty,   \\ 
 \end{array}
\]
provided $f,g$ are absolutely continuous functions defined on [a,b] and $f',g' \in L_\infty  \left[ {a,b} \right]$.

  Recently many researchers have obtained various generalizations, extensions and variants of the above inequalities[2,4-5,8, 10-20,22-24] on time scale calculus. 
German mathematician  Stefan Hilger  has initiated the study of time scales calculus  which unifies the theory of both differential and difference calculus \cite{HIG}. The basic information on time scales can be found in \cite{Agr1,Ana1, Boh1, Boh2, HIG} Motivated by the above literature in this paper we present some dynamic Ostrowski and Cebysev type inequalities on time scales in three variables

    In what follows $\mathbb{R}$ denotes the set of real number and $\mathbb{T}$ is nonempty closed subset of $\mathbb{R}$. For $t \in \mathbb{T}$ the mapping $\sigma ,\rho :\mathbb{T} \to \mathbb{T}$ are defined by $\sigma \left( t \right) = \inf \left\{ {s \in \mathbb{T}:s > t} \right\}$ and $\rho \left( t \right) = \sup \left\{ {s \in \mathbb{T}:s < t} \right\} $ are called the forward and backward jump operators respectively. We define $f:\mathbb{T}\to \mathbb{R}$ is rd-continuous provided f is continuous at each right dense point of $\mathbb{T}$ and has left sided limit at each dense point on $\mathbb{T}$. Now we denote $C_{rd}$ for the set of rd-continuous function defined on $\mathbb{T}$. Let $\mathbb{T}_1,\mathbb{T}_2,\mathbb{T}_3$ be three time scales. Let $\sigma _i ,\rho _i $ and $\Delta _i$ for $i=1,2,3$ denotes the forward jump operator, backward jump operator and delta differentiation operator respectively on $\mathbb{T}_1,\mathbb{T}_2 $ and $\mathbb{T}_3$. Let $a_i<b_i$ be points in $\mathbb{T}_i$ and $[a_i,b_i)$ are half closed bounded intervals in $\mathbb{T}_i$.
		
		We say that a real valued function $f$ on $\mathbb{T}_1  \times \mathbb{T}_2  \times \mathbb{T}_3$ has a $\Delta _1 $ partial derivative $f^{\Delta _1 } \left( {t_1 ,t_2 ,t_3 } \right)$  with respect to $t_1$ if for each $\epsilon >0$ there exists a neighborhood $U_{t_1 }$ of $t_1$ such that 	
\[
\begin{array}{l}
 \left| {f\left( {\sigma _1 (t_1 ),t_2 ,t_3 } \right) - f(s,t_2 ,t_3 ) - f^{\Delta _1 } \left( {t_1 ,t_2 ,t_3 } \right)\left( {\sigma _1 (t_1 ) - s - \tau } \right)} \right| \\ 
  \le \varepsilon \left| {\sigma _1 (t_1 ) - s - \tau } \right|, \\ 
 \end{array}
\]
for each $s \in U_{t_1 }$,	$t_2 \in \mathbb{T}_2$	and $t_3 \in \mathbb{T}_3$.

If $f(x,y,z)$ is delta differential function defined on $\left[ {a_1 ,b_1 } \right] \times \left[ {a_2 ,b_2 } \right] \times \left[ {a_3 ,b_3 } \right]$ for $a_1 ,b_1,a_2 ,b_2,a_3 ,b_3$ in $\mathbb{R}^{+}$. Then its partial delta derivative are defined by \\
$f^{\Delta _1 } \left( {x,y,z} \right) = \frac{\partial }{{\partial x}}f\left( {x,y,z} \right), f^{\Delta _2 } \left( {x,y,z} \right) = \frac{\partial }{{\partial y}}f\left( {x,y,z} \right), f^{\Delta _3 } \left( {x,y,z} \right) = \frac{\partial }{{\partial z}}f\left( {x,y,z} \right)$.
\paragraph{}Let $CC_{rd}$ denotes the set of all rd-continuous functions and let $CC'_{rd}$ denotes the set of all functions for which  $\Delta _1 $ partial derivative, $\Delta _2 $ partial derivative and $\Delta _3 $ partial derivative exists and are in $CC_{rd}$.

\section{\textbf{Ostrowski Inequality In Three Variables on Time Scales}}
Now we give dynamic Ostrowski Inequality in three Variables on Time Scales
\paragraph{\textbf{Theorem 2.1}} Let $f \in CC_{rd} \left( {\left[ {a_1 ,b_1 } \right] \times \left[ {a_2 ,b_2 } \right] \times \left[ {a_3 ,b_3 } \right],\mathbb{R}} \right)$
\begin{align*}
& \left| {\int\limits_{a_1 }^{b_1 } {\int\limits_{a_2 }^{b_2 } {\int\limits_{a_2 }^{b_3 } {f\left( {x,y,z} \right)} } } } \right.\Delta _3 z\Delta _2 y\Delta _1 x \\
&- \frac{1}{8}\left[ {\left( {b_1  - a_1 } \right)\left( {b_2  - a_2 } \right)\left( {b_3  - a_3 } \right)} \right.\left[ {f\left( {\sigma _1 \left( s \right),\sigma _2 \left( t \right),\sigma _3 \left( \tau  \right)} \right)} \right. \\ 
& + f\left( {b_1 ,b_2 ,b_3 } \right) + f\left( {\sigma _1 \left( s \right),\sigma _2 \left( t \right),b_3 } \right) + f\left( {\sigma _1 \left( s \right),b_2 ,\sigma _3 \left( \tau  \right)} \right) \\ 
&+ f\left( {b_1 ,\sigma _2 \left( t \right),\sigma _3 \left( \tau  \right)} \right) + f\left( {\sigma _1 \left( s \right),b_2 ,b_3 } \right) + f\left( {b_1 ,b_2 ,\sigma _3 \left( \tau  \right)} \right) \\
&\left. {\left. { + f\left( {b_1 ,\sigma _2 \left( t \right),b_3 } \right)} \right]} \right] \\ 
&+ \frac{1}{4}\left( {b_2  - a_2 } \right)\left( {b_3  - a_3 } \right)\int\limits_{a_1 }^{b_1 } {\left[ {f\left( {x,\sigma _2 \left( t \right),\sigma _3 \left( \tau  \right)} \right)} \right.}  + f\left( {x,c,d} \right) \\ 
&\left. { + f\left( {x,b_2 ,\sigma _3 \left( \tau  \right)} \right) + f\left( {x,\sigma _2 \left( t \right),b_3 } \right)} \right]\Delta _1 x \\ 
& + \frac{1}{4}\left( {b_1  - a_1 } \right)\left( {b_2  - a_2 } \right)\int\limits_{a_3 }^{b_3 } {\left[ {f\left( {\sigma _1 \left( s \right),\sigma _2 \left( t \right),z} \right)} \right.}  + f\left( {b_1 ,b_2 ,z} \right) \\ 
&\left. { + f\left( {b_1 ,\sigma _2 \left( t \right),z} \right) + f\left( {\sigma _1 \left( s \right),c,z} \right)} \right]\Delta _3 z \\ 
& - \frac{1}{2}\left( {b_1  - a_1 } \right)\int\limits_{a_2 }^{b_2 } {\int\limits_{a_3 }^{b_3 } {\left[ {f\left( {\sigma _1 \left( s \right),y,z} \right) + f\left( {b_1 ,y,z} \right)} \right]} } \Delta _3 z\Delta _2 y \\
& - \frac{1}{2}\left( {b_2  - a_2 } \right)\int\limits_{a_2 }^{b_2 } {\int\limits_{a_3 }^{b_3 } {\left[ {f\left( {x,\sigma _2 \left( t \right),z} \right) + f\left( {x,b_2 ,z} \right)} \right]} } \Delta _3 z\Delta _2 y \\ 
& - \frac{1}{2}\left( {b_3  - a_3 } \right)\int\limits_{a_1 }^{b_1 } {\int\limits_{a_2 }^{b_2 } {\left[ {f\left( {x,y,\sigma _3 \left( \tau  \right)} \right) + f\left( {x,y,b_3 } \right)} \right]} } \Delta _2 y\Delta _1 x\\
& \le \int\limits_{a_1 }^{b_1 } {\int\limits_{a_2 }^{b_2 } {\int\limits_{a_3 }^{b_3 } {\left| {\frac{{\partial ^3 f\left( {\alpha ,\beta ,\gamma } \right)}}{{\Delta _3 \gamma \Delta _2 \beta \Delta _1 \alpha }}} \right|} } } \Delta _3 \gamma \Delta _2 \beta \Delta _1 \alpha. 
\tag{2.1}
\end{align*}
\paragraph{\textbf{Proof.}}:
From the hypotheses we have
\begin{align*}
& \int\limits_{\sigma _1 \left( s \right)}^x {\int\limits_{\sigma _2 \left( t \right)}^y {\int\limits_{\sigma _3 \left( \tau  \right)}^z {\frac{{\partial ^3 f\left( {\alpha ,\beta ,\gamma } \right)}}{{\Delta _3 \gamma \Delta _2 \beta \Delta _1 \alpha }}} } } \Delta _3 \gamma \Delta _2 \beta \Delta _1 \alpha  \\
& = \left. {\int\limits_{\sigma _1 \left( s \right)}^x {\int\limits_{\sigma _2 \left( t \right)}^y {\frac{{\partial ^2 f\left( {\alpha ,\beta ,z} \right)}}{{\Delta _2 \beta \Delta _1 \alpha }}} } } \right|_{\sigma _3 \left( \tau  \right)}^z \Delta _2 \beta \Delta _1 \alpha \\ 
& = \int\limits_{\sigma _1 \left( s \right)}^x {\int\limits_{\sigma _2 \left( t \right)}^y {\frac{{\partial ^2 f\left( {\alpha ,\beta ,z} \right)}}{{\Delta _2 \beta \Delta _1 \alpha }}} } \Delta _2 \beta \Delta _1 \alpha  - \int\limits_{\sigma _1 \left( s \right)}^x {\int\limits_{\sigma _2 \left( t \right)}^y {\frac{{\partial ^2 f\left( {\alpha ,\beta ,\sigma _3 \left( \tau  \right)} \right)}}{{\Delta _2 \beta \Delta _1 \alpha }}} } \Delta _2 \beta \Delta _1 \alpha  \\ 
& = \int\limits_{\sigma _1 \left( s \right)}^x {\left. {\frac{{\partial f\left( {\alpha ,\beta ,z} \right)}}{{\Delta _1 \alpha }}} \right|_{\sigma _2 \left( t \right)}^y } \Delta _1 \alpha  - \int\limits_{\sigma _1 \left( s \right)}^x {\left. {\frac{{\partial f\left( {\alpha ,\beta ,\sigma _3 \left( \tau  \right)} \right)}}{{\Delta _1 \alpha }}} \right|_{\sigma _2 \left( t \right)}^y } \Delta _1 \alpha  \\ 
&  = \int\limits_{\sigma _1 \left( s \right)}^x {\frac{{\partial f\left( {\alpha ,y,z} \right)}}{{\Delta _1 \alpha }}} \Delta _1 \alpha  - \int\limits_{\sigma _1 \left( s \right)}^x {\frac{{\partial f\left( {\alpha ,\sigma _2 \left( t \right),z} \right)}}{{\Delta _1 \alpha }}} \Delta _1 \alpha  \\ 
& - \int\limits_{\sigma _1 \left( s \right)}^x {\frac{{\partial f\left( {\alpha ,y,\sigma _3 \left( \tau  \right)} \right)}}{{\Delta _1 \alpha }}} \Delta _1 \alpha  + \int\limits_{\sigma _1 \left( s \right)}^x {\frac{{\partial f\left( {\alpha ,\sigma _2 \left( t \right),\sigma _3 \left( \tau  \right)} \right)}}{{\Delta _1 \alpha }}} \Delta _1 \alpha  \\ 
& = \left. {f\left( {\alpha ,y,z} \right)} \right|_{\sigma _1 \left( s \right)}^x  - \left. {f\left( {\alpha ,\sigma _2 \left( t \right),z} \right)} \right|_{\sigma _1 \left( s \right)}^x  - \left. {f\left( {\alpha ,y,\sigma _3 \left( \tau  \right)} \right)} \right|_{\sigma _1 \left( s \right)}^x  \\ 
& + \left. {f\left( {\alpha ,\sigma _2 \left( t \right),\sigma _3 \left( t \right)} \right)} \right|_{\sigma _1 \left( s \right)}^x  \\ 
&  = f(x,y,z) - f(\sigma _1 \left( s \right),y,z) - f\left( {x,\sigma _2 \left( t \right),z} \right) \\ 
& + f(\sigma _1 \left( s \right),\sigma _2 \left( t \right),z) - f(x,y,\sigma _3 \left( \tau  \right)) + f\left( {\sigma _1 \left( s \right),y,\sigma _3 \left( \tau  \right)} \right) \\ 
&	+ f\left( {x,\sigma _2 \left( t \right),\sigma _3 \left( \tau  \right)} \right) - f\left( {\sigma _1 \left( s \right),\sigma _2 \left( t \right),\sigma _3 \left( \tau  \right)} \right),
\end{align*}
that is 
\begin{align*}
&f\left( {x,y,z} \right) = f(\sigma _1 \left( s \right),y,z) + f\left( {x,\sigma _2 \left( t \right),z} \right) + f(x,y,\sigma _3 \left( \tau  \right)) \\ 
&+ f\left( {\sigma _1 \left( s \right),\sigma _2 \left( t \right),\sigma _3 \left( \tau  \right)} \right) - f(\sigma _1 \left( s \right),\sigma _2 \left( t \right),z) - f\left( {\sigma _1 \left( s \right),y,\sigma _3 \left( \tau  \right)} \right) \\ 
&- f\left( {x,\sigma _2 \left( t \right),\sigma _3 \left( \tau  \right)} \right) + \int\limits_{\sigma _1 \left( s \right)}^x {\int\limits_{\sigma _2 \left( t \right)}^y {\int\limits_{\sigma _3 \left( \tau  \right)}^z {\frac{{\partial ^3 f\left( {\alpha ,\beta ,\gamma } \right)}}{{\Delta _3 \gamma \Delta _2 \beta \Delta _1 \alpha }}} } } \Delta _3 \gamma \Delta _2 \beta \Delta _1 \alpha. 
\tag{2.2}
\end{align*}
Similarly we have
\begin{align*}
& f\left( {x,y,z} \right) = f(x,y,b_3 ) + f\left( {x,\sigma _2 \left( t \right),z} \right) + f(\sigma _1 \left( s \right),\sigma _2 \left( t \right),b_3 ) \\ 
& - f(\sigma _1 \left( s \right),y,b_3 ) + f(\sigma _1 \left( s \right),y,z) - f(x,\sigma _2 \left( t \right),b_3 ) \\ 
&- f\left( {\sigma _1 \left( s \right),\sigma _2 \left( t \right),z} \right) - \int\limits_{\sigma _1 \left( s \right)}^x {\int\limits_{\sigma _2 \left( t \right)}^y {\int\limits_z^{b_3 } {\frac{{\partial ^3 f\left( {\alpha ,\beta ,\gamma } \right)}}{{\Delta _3 \gamma \Delta _2 \beta \Delta _1 \alpha }}} } } \Delta _3 \gamma \Delta _2 \beta \Delta _1 \alpha,  
\tag{2.3}
\end{align*}
\begin{align*}
&f\left( {x,y,z} \right) = f\left( {x,b_2 ,z} \right) + f\left( {\sigma _1 \left( s \right),b_2 ,\sigma _3 \left( \tau  \right)} \right) + f\left( {x,y,\sigma _3 \left( \tau  \right)} \right) \\ 
& + f\left( {\sigma _1 \left( s \right),y,z} \right) - f\left( {\sigma _1 \left( s \right),b_2 ,z} \right) - f\left( {x,b_2 ,\sigma _3 \left( \tau  \right)} \right) \\ 
&- f\left( {\sigma _1 \left( s \right),y,\sigma _3 \left( \tau  \right)} \right) - \int\limits_{\sigma _1 \left( s \right)}^x {\int\limits_y^{b_2 } {\int\limits_{\sigma _3 \left( \tau  \right)}^z {\frac{{\partial ^3 f\left( {\alpha ,\beta ,\gamma } \right)}}{{\Delta _3 \gamma \Delta _2 \beta \Delta _1 \alpha }}} } } \Delta _3 \gamma \Delta _2 \beta \Delta _1 \alpha  ,
\tag{2.4}
\end{align*}
\begin{align*}
&f\left( {x,y,z} \right) = f\left( {b_1 ,y,z} \right) + f\left( {x,\sigma _2 \left( t \right),z} \right) + f\left( {x,y,\sigma _3 \left( \tau  \right)} \right) \\ 
&+ f\left( {b_1 ,\sigma _2 \left( t \right),\sigma _3 \left( \tau  \right)} \right) - f\left( {b_1 ,\sigma _2 \left( t \right),z} \right) - f\left( {b,y,\sigma _3 \left( \tau  \right)} \right) \\ 
& - f\left( {x,\sigma _2 \left( t \right),\sigma _3 \left( \tau  \right)} \right) - \int\limits_x^{b_1 } {\int\limits_{\sigma _2 \left( t \right)}^y {\int\limits_{\sigma _3 \left( \tau  \right)}^z {\frac{{\partial ^3 f\left( {\alpha ,\beta ,\gamma } \right)}}{{\Delta _3 \gamma \Delta _2 \beta \Delta _1 \alpha }}} } } \Delta _3 \gamma \Delta _2 \beta \Delta _1 \alpha  ,
\tag{2.5}
\end{align*}
\begin{align*}
&f\left( {x,y,z} \right) = f\left( {\sigma _1 \left( s \right),b_2 ,b_3 } \right) + f\left( {x,y,b_3 } \right) + f\left( {x,b_2 ,z} \right) \\ 
& + f\left( {\sigma _1 \left( s \right),y,z} \right) - f\left( {x,b_2 ,b_3 } \right) - f\left( {\sigma _1 \left( s \right),y,b_3 } \right) \\
& - f\left( {\sigma _1 \left( s \right),b_2 ,z} \right) + \int\limits_{\sigma _1 \left( s \right)}^{x } {\int\limits_y^{b_2 } {\int\limits_z^{b_3 } {\frac{{\partial ^3 f\left( {\alpha ,\beta ,\gamma } \right)}}{{\Delta _3 \gamma \Delta _2 \beta \Delta _1 \alpha }}} } } \Delta _3 \gamma \Delta _2 \beta \Delta _1 \alpha,   
\tag{2.6}
\end{align*}
\begin{align*}
&f\left( {x,y,z} \right) = f\left( {x,b_2 ,z} \right) + f\left( {b_1 ,y,z} \right) + f\left( {b_1 ,b_2 ,\sigma _3 \left( \tau  \right)} \right) \\
&+ f\left( {x,y,\sigma _3 \left( \tau  \right)} \right) - f\left( {b_1 ,b_2 ,z} \right) - f\left( {x,b_2 ,\sigma _3 \left( \tau  \right)} \right) \\
&- f\left( {b_1 ,y,\sigma _3 \left( \tau  \right)} \right) + \int\limits_x^{b_1 } {\int\limits_y^{b_2 } {\int\limits_{\sigma _3 \left( \tau  \right)}^z {\frac{{\partial ^3 f\left( {\alpha ,\beta ,\gamma } \right)}}{{\Delta _3 \gamma \Delta _2 \beta \Delta _1 \alpha }}} } } \Delta _3 \gamma \Delta _2 \beta \Delta _1 \alpha,  \\ 
\tag{2.7}
\end{align*}
\begin{align*}
&f\left( {x,y,z} \right) = f\left( {x,y,b_3 } \right) + f\left( {b_1 ,\sigma _2 \left( t \right),b_3 } \right) + f\left( {b_1 ,y,z} \right) \\
&+ f\left( {x,\sigma _2 \left( t \right),z} \right) - f\left( {b_1 ,y,b_3 } \right) - f\left( {x,\sigma _2 \left( t \right),b_3 } \right) \\
& - f\left( {b_1 ,\sigma _2 \left( t \right),y} \right) + \int\limits_x^{b_1 } {\int\limits_{\sigma _2 \left( t \right)}^y {\int\limits_z^{b_3 } {\frac{{\partial ^3 f\left( {\alpha ,\beta ,\gamma } \right)}}{{\Delta _3 \gamma \Delta _2 \beta \Delta _1 \alpha }}\Delta _3 \gamma \Delta _2 \beta \Delta _1 \alpha } } } , \\ 
\tag{2.8}
\end{align*}
\begin{align*}
&f\left( {x,y,z} \right) = f\left( {b_1 ,b_2 ,b_3 } \right) + f\left( {x,y,b_3 } \right) + f\left( {x,b_2 ,z} \right) \\ 
&+ f\left( {b_1 ,y,z} \right) + f\left( {x,b_2 ,b_3 } \right) - f\left( {b_1 ,y,b_3 } \right) \\ 
&- f\left( {x,b_2 ,z} \right) - \int\limits_x^{b_1 } {\int\limits_y^{b_2 } {\int\limits_z^{b_3 } {\frac{{\partial ^3 f\left( {\alpha ,\beta ,\gamma } \right)}}{{\Delta _3 \gamma \Delta _2 \beta \Delta _1 \alpha }}\Delta _3 \gamma \Delta _2 \beta \Delta _1 \alpha } } } .
\tag{2.9}
\end{align*}
Integrating the above equations(2.2-2.9) and adding them we get the required inequality (2.1).
\paragraph{}If we have $\mathbb{T}=\mathbb{R}$ in above theorem we get the continuous version to the above inequality. 
\paragraph{\textbf{Corollary 2.2}} Let $f:\left[ {a_1 ,b_1 } \right] \times \left[ {a_2 ,b_2 } \right] \times \left[ {a_3 ,b_3 } \right] \to \mathbb{R}$ be differentiable and continuous function. Then
\begin{align*}
& \left| {\int\limits_{a_1 }^{b_1 } {\int\limits_{a_2 }^{b_2 } {\int\limits_{a_3 }^{b_3 } {f\left( {x,y,z} \right)dzdydx} } } } \right. \\ 
& - \frac{1}{8}\left[ {\left( {b_1  - a_1 } \right)} \right.\left( {b_2  - a_2 } \right)\left( {b_3  - a_3 } \right)\left[ {f\left( {s,t,\tau } \right)} \right. + f\left( {b_1 ,b_2 ,b_3 } \right) \\ 
&+ f\left( {s,t,b_3 } \right) + f\left( {s,b_2 ,\tau } \right) + f\left( {b_1 ,t,\tau } \right) \\ 
&\left. {\left. { + f\left( {s,b_2 ,b_3 } \right) + f\left( {b_1 ,b_2 ,\tau } \right) + f\left( {b_1 ,t,b_3 } \right)} \right]} \right] \\ 
&+ \frac{1}{4}\left( {b_2  - a_2 } \right)\left( {b_3  - a_3 } \right)\int\limits_{a_1 }^{b_1 } {\left[ {f\left( {x,t,\tau } \right)} \right.}  + f\left( {x,c,d} \right) \\ 
& \left. { + f\left( {x,b_2 ,\tau } \right) + f\left( {x,b_2 ,b_3 } \right)} \right]dx \\ 
&+ \frac{1}{4}\left( {b_1  - a_1 } \right)\left( {b_2  - a_2 } \right)\int\limits_{a_1 }^{b_1 } {\left[ {f\left( {s,t,z} \right)} \right.}  + f\left( {b_1 ,b_2 ,z} \right) \\ 
&\left. { + f\left( {b_1 ,t,z} \right) + f\left( {s,c,z} \right)} \right]dz \\ 
&- \frac{1}{2}\left( {b_1  - a_1 } \right)\int\limits_{a_2 }^{b_2 } {\int\limits_{a_3 }^{b_3 } {\left[ {f\left( {s,y,z} \right) + f\left( {b_1 ,y,z} \right)} \right]} } dzdy \\ 
&- \frac{1}{2}\left( {b_2  - a_2 } \right)\int\limits_{a_1 }^{b_1 } {\int\limits_{a_3 }^{b_3 } {\left[ {f\left( {x,t,z} \right) + f\left( {x,b_2 ,z} \right)} \right]} } dzdx \\ 
& - \frac{1}{2}\left( {b_3  - a_3 } \right)\int\limits_{a_1 }^{b_1 } {\int\limits_{a_2 }^{b_2 } {\left[ {f\left( {x,y,\tau } \right) + f\left( {x,y,b_3 } \right)} \right]} } dydx |\\ 
&\le \int\limits_{a_1 }^{b_1 } {\int\limits_{a_2 }^{b_2 } {\int\limits_{a_3 }^{b_3 } {\left| {D_3 D_2 D_1 f\left( {\alpha ,\beta ,\gamma } \right)} \right|d\gamma d\beta d\alpha } } }
\tag{2.10}
\end{align*}
\paragraph{}If we have $\mathbb{T}=\mathbb{Z}$ in above theorem we get the discrete version of above inequality
\paragraph{\textbf{Corollary 2.3}} Let $ \mathbb{Z}_1  = \left\{ {1,2,...,b_1  + 1} \right\}$,$\mathbb{Z}_2  = \left\{ {1,2,...,b_2  + 1} \right\}$,$\mathbb{Z}_3  = \left\{ {1,2,...,b_3  + 1} \right\}$ for $b_1,b_2,b_3 \in \mathbb{Z}$ and $f:\mathbb{Z}_1  \times \mathbb{Z}_2  \times \mathbb{Z}_3  \to \mathbb{R}$
\begin{align*}
&\left| {\sum\limits_{\alpha  = 1}^{b_1 } {\sum\limits_{\beta  = 1}^{b_2 } {\sum\limits_{\gamma  = 1}^{b_3 } {f\left( {\alpha ,\beta ,\gamma } \right)} } } } \right. - \frac{1}{8}\left[ {f\left( {1,1,1} \right)} \right. + f\left( {b_1  + 1,b_2  + 1,b_3  + 1} \right) \\ 
& + f\left( {1,1,b_3  + 1} \right) + f\left( {1,b_2  + 1,1} \right) + f\left( {b_1  + 1,1,1} \right) \\ 
&\left. { + f\left( {1,b_2  + 1,b_3  + 1} \right) + f\left( {b_1  + 1,b_2  + 1,1} \right) + f\left( {b_1  + 1,1,b_3  + 1} \right)} \right] \\ 
&+ \frac{1}{4}b_2 b_3 \sum\limits_{r = 1}^{b_1 } {\left[ {f\left( {r,1,1} \right)} \right.}  + f\left( {r,1,b_3  + 1} \right) + f\left( {r,b_2  + 1,1} \right) \\ 
&\left. { + f\left( {r,b_2  + 1,b_3  + 1} \right)} \right] \\ 
&+ \frac{1}{4}b_1 b_3 \sum\limits_{r = 1}^{b_2 } {\left[ {f\left( {b_1  + 1,1,b_3  + 1} \right)} \right.}  + f\left( {b_1  + 1,s,1} \right) + f\left( {1,s,b_3  + 1} \right) \\ 
&\left. { + f\left( {r,b_2  + 1,b_3  + 1} \right)} \right] \\
&  + \frac{1}{4}b_1 b_2 \sum\limits_{r = 1}^{b_2 } {\left[ {f\left( {b_1  + 1,b_2  + 1,t} \right)} \right.}  + f\left( {b_1  + 1,1,t} \right) + f\left( {1,b_2  + 1,t} \right) \\
&\left. { + f\left( {1,1,t} \right)} \right]  \\
&- \frac{1}{2}b_1 \sum\limits_{s = 1}^{b_2 } {\sum\limits_{t = 1}^{b_3 } {\left[ {f\left( {1,s,t} \right) + f\left( {b_1  + 1,s,t} \right)} \right]} }  \\
&- \frac{1}{2}b_2 \sum\limits_{r = 1}^{b_1 } {\sum\limits_{s = 1}^{b_2 } {\left[ {f\left( {r,1,t} \right) + f\left( {r,b_2  + 1,t} \right)} \right]} }  \\ 
&- \frac{1}{2}b_3 \sum\limits_{r = 1}^{b_1 } {\sum\limits_{s = 1}^{b_2 } {\left[ {f\left( {r,s,1} \right) + f\left( {r,s,b_3  + 1} \right)} \right]} }  \\ 
& \le \frac{1}{8}b_1 b_2 b_3 \sum\limits_{\alpha  = 1}^{b_1 } {\sum\limits_{\beta  = 1}^{b_2 } {\sum\limits_{\gamma  = 1}^{b_3 } {\left| {\Delta _3 \Delta _2 \Delta _1 f\left( {\alpha ,\beta ,\gamma } \right)} \right|} } }  \\ 
\tag{2.11}
\end{align*}
\section{\textbf{Cebysev type inequalities in three Variables}}
Now we give the Cebysev type inequalities in three Variables on time scales. 
\paragraph{}We use following notations to simplify the details.
\begin{align*}
A\left( {f\left( {x,y,z} \right)} \right)
&= \frac{1}{8}\left[ {f\left( {\sigma _1 \left( s \right),\sigma _2 \left( t \right),\sigma _3 \left( \tau  \right)} \right) + f\left( {b_1 ,b_2 ,b_3 } \right)} \right. \\
& + f\left( {\sigma _1 \left( s \right),\sigma _2 \left( t \right),b_3 } \right) + f\left( {\sigma _1 \left( s \right),b_2 ,\sigma _3 \left( \tau  \right)} \right) \\ 
&+ f\left( {b_1 ,\sigma _2 \left( t \right),\sigma _3 \left( \tau  \right)} \right) + f\left( {\sigma _1 \left( s \right),b_2 ,b_3 } \right) \\ 
&\left. { + f\left( {b_1 ,b_2 ,\sigma _3 \left( \tau  \right)} \right) + f\left( {b_1 ,\sigma _2 \left( t \right),b_3 } \right)} \right] \\ 
& - \frac{1}{4}\left[ {f\left( {x,\sigma _2 \left( t \right),\sigma _3 \left( \tau  \right)} \right)} \right. + f\left( {x,b_2 ,b_3 } \right) \\
&\left. { + f\left( {x,b_2 ,\sigma _3 \left( \tau  \right)} \right) + f\left( {x,\sigma _2 \left( t \right),b_3 } \right)} \right] \\ 
& - \frac{1}{4}\left[ {f\left( {\sigma _1 \left( s \right),y,\sigma _3 \left( \tau  \right)} \right) + f\left( {b_1 ,y,b_3 } \right)} \right. \\ 
& \left. { + f\left( {\sigma _1 \left( s \right),y,b_3 } \right) + f\left( {b_1 ,y,\sigma _3 \left( \tau  \right)} \right)} \right] \\ 
& - \frac{1}{4}\left[ {f\left( {\sigma _1 \left( s \right),\sigma _2 \left( t \right),z} \right) + f\left( {b_1 ,c,z} \right)} \right. \\ 
&\left. { + f\left( {b_1 ,\sigma _2 \left( t \right),z} \right) + f\left( {\sigma _1 \left( s \right),b_2 ,z} \right)} \right] \\ 
&+ \frac{1}{2}\left[ {f\left( {\sigma _1 \left( s \right),y,z} \right) + f\left( {b_1 ,y,z} \right)} \right] \\ 
& + \frac{1}{2}\left[ {f\left( {x,\sigma _2 \left( t \right),z} \right) + f\left( {x,b_2 ,z} \right)} \right] \\ 
& + \frac{1}{2}\left[ {f\left( {x,y,\sigma _3 \left( \tau  \right)} \right) + f\left( {x,y,b_3 } \right)} \right], 
\tag{3.1}
\end{align*}
and
\begin{align*}
B\left( {f^{\Delta _3 \Delta _2 \Delta _1 } \left( {x,y,z} \right)} \right)  
&  = \int\limits_{\sigma _1 \left( s \right)}^x {\int\limits_{\sigma _2 \left( t \right)}^y {\int\limits_{\sigma _3 \left( \tau  \right)}^z {\frac{{\partial ^3 f\left( {\alpha ,\beta ,\gamma } \right)}}{{\Delta _3 \gamma \Delta _2 \beta \Delta _1 \alpha }}} } } \Delta _3 \gamma \Delta _2 \beta \Delta _1 \alpha  \\ 
&- \int\limits_{\sigma _1 \left( s \right)}^x {\int\limits_{\sigma _2 \left( t \right)}^y {\int\limits_z^{b_3 } {\frac{{\partial ^3 f\left( {\alpha ,\beta ,\gamma } \right)}}{{\Delta _3 \gamma \Delta _2 \beta \Delta _1 \alpha }}} } } \Delta _3 \gamma \Delta _2 \beta \Delta _1 \alpha  \\ 
& - \int\limits_{\sigma _1 \left( s \right)}^x {\int\limits_y^{b_2 } {\int\limits_{\sigma _3 \left( \tau  \right)}^z {\frac{{\partial ^3 f\left( {\alpha ,\beta ,\gamma } \right)}}{{\Delta _3 \gamma \Delta _2 \beta \Delta _1 \alpha }}} } } \Delta _3 \gamma \Delta _2 \beta \Delta _1 \alpha  \\
&- \int\limits_x^{b_1 } {\int\limits_{\sigma _2 \left( t \right)}^y {\int\limits_{\sigma _3 \left( \tau  \right)}^z {\frac{{\partial ^3 f\left( {\alpha ,\beta ,\gamma } \right)}}{{\Delta _3 \gamma \Delta _2 \beta \Delta _1 \alpha }}} } } \Delta _3 \gamma \Delta _2 \beta \Delta _1 \alpha  \\ 
& + \int\limits_{\sigma _1 \left( s \right)}^x {\int\limits_y^{b_2 } {\int\limits_z^{b_3 } {\frac{{\partial ^3 f\left( {\alpha ,\beta ,\gamma } \right)}}{{\Delta _3 \gamma \Delta _2 \beta \Delta _1 \alpha }}} } } \Delta _3 \gamma \Delta _2 \beta \Delta _1 \alpha  \\ 
&+ \int\limits_x^{b_1 } {\int\limits_y^{b_2 } {\int\limits_{\sigma _3 \left( \tau  \right)}^z {\frac{{\partial ^3 f\left( {\alpha ,\beta ,\gamma } \right)}}{{\Delta _3 \gamma \Delta _2 \beta \Delta _1 \alpha }}} } } \Delta _3 \gamma \Delta _2 \beta \Delta _1 \alpha  \\ 
&+ \int\limits_x^{b_1 } {\int\limits_{\sigma _2 \left( t \right)}^y {\int\limits_z^{b_3 } {\frac{{\partial ^3 f\left( {\alpha ,\beta ,\gamma } \right)}}{{\Delta _3 \gamma \Delta _2 \beta \Delta _1 \alpha }}} } } \Delta _3 \gamma \Delta _2 \beta \Delta _1 \alpha  \\ 
&- \int\limits_x^{b_1 } {\int\limits_y^{b_2 } {\int\limits_z^{b_3 } {\frac{{\partial ^3 f\left( {\alpha ,\beta ,\gamma } \right)}}{{\Delta _3 \gamma \Delta _2 \beta \Delta _1 \alpha }}} } } \Delta _3 \gamma \Delta _2 \beta \Delta _1 \alpha.  
\tag{3.2}
\end{align*}
\paragraph{}Now we give Cebysev type inequality on time scales involving functions in three variables.
 \paragraph{\textbf{Theorem 3.1}}
Let $f,g \in CC_{rd} \left( {\left[ {a_1 ,b_1 } \right] \times \left[ {a_2 ,b_2 } \right] \times \left[ {a_3 ,b_3 } \right],\mathbb{R}} \right)$ be rd-continuous function and $f^{\Delta _3 \Delta _2 \Delta _1 } \left( {x,y,z} \right)$ and $g^{\Delta _3 \Delta _2 \Delta _1 } \left( {x,y,z} \right)$ exist and rd-continuous and bounded. Then
\begin{align*}
&\int\limits_{\sigma _1 \left( s \right)}^{b_1 } {\int\limits_{\sigma _2 \left( t \right)}^{b_2 } {\int\limits_{\sigma _3 \left( \tau  \right)}^{b_3 } {\left[ {f\left( {x,y,z} \right)g\left( {x,y,z} \right)} \right.} } }  - \frac{1}{2}\left[ {f\left( {x,y,z} \right)A\left( {g\left( {x,y,z} \right)} \right)} \right. \\ 
&\left. { + g\left( {x,y,z} \right)A\left( {f\left( {x,y,z} \right)} \right)} \right]\Delta _3 z\Delta _2 y\Delta _1 x \\ 
&\le \frac{1}{{16}}\left( {b_1  - \sigma _1 \left( s \right)} \right)\left( {b_2  - \sigma _2 \left( t \right)} \right)\left( {b_3  - \sigma _3 \left( \tau  \right)} \right) \\ 
&\int\limits_{\sigma _1 \left( s \right)}^{b_1 } {\int\limits_{\sigma _2 \left( t \right)}^{b_2 } {\int\limits_{\sigma _3 \left( \tau  \right)}^{b_3 } {\left[ {\left| {g\left( {x,y,z} \right)} \right|\left\| {f^{\Delta _3 \Delta _2 \Delta _1 } } \right\|_\infty  } \right.} } }  \\
&\left. { + \left| {f\left( {x,y,z} \right)} \right|\left\| {g^{\Delta _3 \Delta _2 \Delta _1 } } \right\|_\infty  } \right]\Delta _3 z\Delta _2 y\Delta _1 x.
\tag{3.3}
\end{align*}
\paragraph{\textbf{Proof.}}:
Adding $(2.2)-(2.9)$ and using the notations in equation $(3.1)$ and $(3.2)$ we have 
\[
f\left( {x,y,z} \right) - A\left( {f\left( {x,y,z} \right)} \right) = \frac{1}{8}B\left( {f^{\Delta _3 \Delta _2 \Delta _1 } \left( {x,y,z} \right)} \right),
\tag{3.4}
\]
for $(x,y,z) \in [a_1 ,b_1] \times [a_2 ,b_2] \times [a_3 ,b_3]$. \\
Similarly for $g$ we have
\[
g\left( {x,y,z} \right) - A\left( {g\left( {x,y,z} \right)} \right) = \frac{1}{8}B\left( {g^{\Delta _3 \Delta _2 \Delta _1 } \left( {x,y,z} \right)} \right).
\tag{3.5}
\]
Multiplying $(3.4)$ by $g(x,y,z)$ and $(3.5)$ by $f(x,y,z)$ and adding resultant identities and then integrating we have
\begin{align*}
& \int\limits_{\sigma _1 \left( s \right)}^{b_1 } {\int\limits_{\sigma _2 \left( t \right)}^{b_2 } {\int\limits_{\sigma _3 \left( \tau  \right)}^{b_3 } {\left[ {f\left( {x,y,z} \right)g\left( {x,y,z} \right)} \right.} } }  - \frac{1}{2}\left[ {f\left( {x,y,z} \right)A\left( {g\left( {x,y,z} \right)} \right)} \right. \\
&\left. { + g\left( {x,y,z} \right)A\left( {f\left( {x,y,z} \right)} \right)} \right]\Delta _3 z\Delta _2 y\Delta _1 x \\ 
&= \frac{1}{{16}}\int\limits_{\sigma _1 \left( s \right)}^{b_1 } {\int\limits_{\sigma _2 \left( t \right)}^{b_2 } {\int\limits_{\sigma _3 \left( \tau  \right)}^{b_3 } {\left[ {g\left( {x,y,z} \right)B\left( {f^{\Delta _3 \Delta _2 \Delta _1 } \left( {x,y,z} \right)} \right)} \right.} } }  \\ 
&\left. { + f\left( {x,y,z} \right)B\left( {g^{\Delta _3 \Delta _2 \Delta _1 } \left( {x,y,z} \right)} \right)} \right]\Delta _3 z\Delta _2 y\Delta _1 x.
\tag{3.6}
\end{align*}
From the properties of modulus and integrals we have
\begin{align*}
&\left| {B\left( {f^{\Delta _3 \Delta _2 \Delta _1 } \left( {x,y,z} \right)} \right)} \right| \\
&\le \int\limits_{\sigma _1 \left( s \right)}^{b_1 } {\int\limits_{\sigma _2 \left( t \right)}^{b_2 } {\int\limits_{\sigma _3 \left( \tau  \right)}^{b_3 } {\left| {f^{\Delta _3 \Delta _2 \Delta _1 } \left( {x,y,z} \right)} \right|} } } \Delta _3 z\Delta _2 y\Delta _1 x \\ 
&\le \left\| {f^{\Delta _3 \Delta _2 \Delta _1 } } \right\|_\infty  \left( {b_1  - \sigma _1 \left( s \right)} \right)\left( {b_2  - \sigma _2 \left( t \right)} \right)\left( {b_3  - \sigma _3 \left( \tau  \right)} \right) .
\tag{3.7}
\end{align*}
Similarly we have 
\begin{align*}
&\left| {B\left( {g^{\Delta _3 \Delta _2 \Delta _1 } \left( {x,y,z} \right)} \right)} \right| \\ 
&\le \left\| {g^{\Delta _3 \Delta _2 \Delta _1 } } \right\|_\infty  \left( {b_1  - \sigma _1 \left( s \right)} \right)\left( {b_2  - \sigma _2 \left( t \right)} \right)\left( {b_3  - \sigma _3 \left( \tau  \right)} \right).
\tag{3.8}
\end{align*}
Now from $(3.6), (3.7),(3.8)$ we have
\begin{align*}
&\int\limits_{\sigma _1 \left( s \right)}^{b_1 } {\int\limits_{\sigma _2 \left( t \right)}^{b_2 } {\int\limits_{\sigma _3 \left( \tau  \right)}^{b_3 } {\left[ {f\left( {x,y,z} \right)g\left( {x,y,z} \right)} \right.} } }  - \frac{1}{2}\left[ {f\left( {x,y,z} \right)A\left( {g\left( {x,y,z} \right)} \right)} \right. \\ 
&\left. { + g\left( {x,y,z} \right)A\left( {f\left( {x,y,z} \right)} \right)} \right]\Delta _3 z\Delta _2 y\Delta _1 x \\ 
&\le \int\limits_{\sigma _1 \left( s \right)}^{b_1 } {\int\limits_{\sigma _2 \left( t \right)}^{b_2 } {\int\limits_{\sigma _3 \left( \tau  \right)}^{b_3 } {\left[ {\left| {g\left( {x,y,z} \right)} \right|} \right.} } } \int\limits_{\sigma _1 \left( s \right)}^{b_1 } {\int\limits_{\sigma _2 \left( t \right)}^{b_2 } {\int\limits_{\sigma _3 \left( \tau  \right)}^{b_3 } {\left| {f^{\Delta _3 \Delta _2 \Delta _1 } \left( {\xi ,\eta ,\gamma } \right)} \right|} } } \Delta _3 \gamma \Delta _2 \eta \Delta _1 \xi  \\ 
&\left. { + \left| {f\left( {x,y,z} \right)} \right|\int\limits_{\sigma _1 \left( s \right)}^{b_1 } {\int\limits_{\sigma _2 \left( t \right)}^{b_2 } {\int\limits_{\sigma _3 \left( \tau  \right)}^{b_3 } {\left| {g^{\Delta _3 \Delta _2 \Delta _1 } \left( {\xi ,\eta ,\gamma } \right)} \right|} } } \Delta _3 \gamma \Delta _2 \eta \Delta _1 \xi } \right]\Delta _3 z\Delta _2 y\Delta _1 x \\ 
&\le \frac{1}{{16}}\left( {b_1  - \sigma _1 \left( s \right)} \right)\left( {b_2  - \sigma _2 \left( t \right)} \right)\left( {b_3  - \sigma _3 \left( \tau  \right)} \right) \\ 
&\int\limits_{\sigma _1 \left( s \right)}^{b_1 } {\int\limits_{\sigma _2 \left( t \right)}^{b_2 } {\int\limits_{\sigma _3 \left( \tau  \right)}^{b_3 } {\left| {g\left( {x,y,z} \right)} \right|\left\| {f^{\Delta _3 \Delta _2 \Delta _1 } } \right\|_\infty  } } }  \\ 
&+\left. {f\left( {x,y,z} \right)\left\| {g^{\Delta _3 \Delta _2 \Delta _1 } } \right\|_\infty  } \right]\Delta _3 z\Delta _2 y\Delta _1 x,
\tag{3.9}
\end{align*}
which is required inequality $(3.2)$.
If we have $\mathbb{T}=\mathbb{R}$ in above theorem we have
\paragraph{\textbf{Corollary 3.2}} Let $f,g:\left[ {a_1 ,b_1 } \right] \times \left[ {a_2 ,b_2 } \right] \times \left[ {a_3 ,b_3 } \right] \to \mathbb{R}$ be differentiable and continuous function and $D_3 D_2 D_1 f\left( {x,y,z} \right)$, 
$D_3 D_2 D_1 g\left( {x,y,z} \right)$ exist continuous and bounded. Then
\begin{align*}
&\left| {\int\limits_a^k {\int\limits_b^m {\int\limits_c^n {\left[ {f\left( {x,y,z} \right)} \right.} } } g\left( {x,y,z} \right) - \left\{ {f\left( {x,y,z} \right)P\left( {g\left( {x,y,z} \right)} \right)} \right.} \right. \\ 
&\left. {\left. {\left. { + g\left( {x,y,z} \right)P\left( {f\left( {x,y,z} \right)} \right)} \right\}} \right]dzdydx} \right| \\ 
&\le \frac{1}{{16}}\left( {k - a} \right)\left( {m - b} \right)\left( {n - c} \right)\int\limits_a^k {\int\limits_b^m {\int\limits_c^n {\left[ {\left| {g\left( {x,y,z} \right)} \right|} \right.} } } \left\| {D_3 D_2 D_1 f} \right\|_\infty   \\ 
& \left. { + \left| {f\left( {x,y,z} \right)} \right|\left\| {D_3 D_2 D_1 g} \right\|_\infty  } \right]dzdydx,
\tag{3.10}
\end{align*}
where 
\begin{align*}
&P\left( {f\left( {x,y,z} \right)} \right) \\ 
&= \frac{1}{8}\left[ {f\left( {a,b,c} \right) + f\left( {k,m,n} \right)} \right. + f(a,m,n) + f(a,b,n) \\ 
&+ f(k,b,n) + f(k,m,c) + f(k,b,c)\left. { + f(a,m,c)} \right] \\ 
& - \frac{1}{4}\left[ {f\left( {x,b,c} \right) + f\left( {x,m,n} \right) + f\left( {x,m,c} \right) + f\left( {x,b,n} \right)} \right] \\
&- \frac{1}{4}\left[ {f\left( {a,y,c} \right) + f\left( {k,y,n} \right) + f\left( {a,y,n} \right) + f\left( {k,y,c} \right)} \right] \\ 
&- \frac{1}{4}\left[ {f\left( {a,b,z} \right) + f\left( {k,m,z} \right) + f\left( {k,b,z} \right) + f\left( {a,m,z} \right)} \right] \\ 
& + \frac{1}{2}\left[ {f\left( {a,y,z} \right) + f\left( {k,y,z} \right)} \right] + \frac{1}{2}\left[ {f\left( {x,b,z} \right) + f\left( {x,m,z} \right)} \right] \\ 
&+ \frac{1}{2}\left[ {f\left( {x,y,c} \right) + f\left( {x,y,n} \right)} \right]. 
\tag{3.10}
\end{align*}
If we have $\mathbb{T}=\mathbb{Z}$ in above theorem we have discrete version of above inequality.
\paragraph{\textbf{Corollary 3.3}}
Let $ Z_1  = \left\{ {1,2,...,k  + 1} \right\}$,$\mathbb{Z}_2  = \left\{ {1,2,...,m  + 1} \right\}$,$\mathbb{Z}_3  = \left\{ {1,2,...,n  + 1} \right\}$ for $b_1,b_2,b_3 \in \mathbb{Z}$ and $f,g:\mathbb{Z}_1  \times \mathbb{Z}_2  \times \mathbb{Z}_3  \to \mathbb{R}$
be functions such that $\Delta _3 \Delta _2 \Delta _1 f\left( {x,y,z} \right)$, $\Delta _3 \Delta _2 \Delta _1 g\left( {x,y,z} \right)$ exist and are bounded. Then
\begin{align*}
&\left| {\sum\limits_{x = 1}^k {\sum\limits_{y = 1}^m {\sum\limits_{z = 1}^m {\left[ {f\left( {x,y,z} \right)g\left( {x,y,z} \right)} \right.}  - \left\{ {f\left( {x,y,z} \right)} \right.} } } \right.Q\left( {g\left( {x,y,z} \right)} \right) \\ 
&\left. {\left. {\left. { + g\left( {x,y,z} \right)Q\left( {f\left( {x,y,z} \right)} \right)} \right]} \right\}} \right| \\ 
& \le \frac{1}{{16}}kmn\sum\limits_{x = 1}^k {\sum\limits_{y = 1}^m {\sum\limits_{z = 1}^m {\left[ {\left| {f\left( {x,y,z} \right)} \right|\left\| {\Delta _3 \Delta _2 \Delta _1 g} \right\|_\infty  } \right.} } }  \\ 
& \left. { + \left| {g\left( {x,y,z} \right)} \right|\left\| {\Delta _3 \Delta _2 \Delta _1 f} \right\|_\infty  } \right],
\tag{3.11}
\end{align*}
where
\begin{align*}
&Q(f(x,y,z)) = \frac{1}{8}\left[ {f\left( {1,1,1} \right) + \left. {e\left( {k + 1,m + 1,n + 1} \right)} \right]} \right. \\
& - \frac{1}{4}\left[ {f\left( {x,1,1} \right) + f(x,1,n + 1) + f(x,m + 1,1) + f(x,m + 1,n + 1)} \right] \\ 
& - \frac{1}{4}\left[ {f\left( {k + 1,y,n + 1} \right) + f(k + 1,y,1) + f(1,y,n + 1) + f(1,y,1)} \right] \\ 
&- \frac{1}{4}\left[ {f\left( {k + 1,m + 1,z} \right) + f(k + 1,z) + f(1,m + 1,t) + f(1,1,t)} \right] \\ 
&+ \frac{1}{2}\left[ {f\left( {1,y,z} \right) + f(k + 1,y,z)} \right] + \frac{1}{2}\left[ {f\left( {x,1,z} \right) + f(x,m + 1,z)} \right] \\ 
& + \frac{1}{2}\left[ {f\left( {x,y,1} \right) + f(x,y,n + 1)} \right].\\
\tag{3.12}
\end{align*}
 \paragraph{\textbf{Acknowledgment}}
This research is supported by Science and Engineering Research Board (SERB, New
Delhi, India) Project File No. SB/S4/MS:861/13.
  
\end{document}